\documentclass[]{article}

\usepackage{amsmath}

\usepackage{amssymb}
\usepackage{amsfonts}
\usepackage{amsthm}
\usepackage{mathtools}

\newtheorem{theorem}{Theorem}[section]
\newtheorem{proposition}[theorem]{Proposition}
\newtheorem{lemma}[theorem]{Lemma}

\theoremstyle{remark}
\newtheorem*{remark}{Remark}
\theoremstyle{definition}

\numberwithin{equation}{section}

\DeclarePairedDelimiter\abs{\lvert}{\rvert}
\DeclarePairedDelimiter\norm{\lVert}{\rVert}

\newcommand{\E}{\mathbb{E}}
\renewcommand{\H}{\mathbb{H}}

\renewcommand{\P}{\mathbb{P}}

\newcommand{\R}{\mathbb{R}}
\newcommand{\N}{\mathbb{N}}

\newcommand{\C}{\mathbb{C}}

\newcommand{\I}{\mathbb{I}}

\newcommand{\normHR}[2]{\norm{#2}_{\dot{H}^{#1}(\R)}}

\renewcommand{\H}{\mathbb{H}}
\DeclareMathOperator{\Conf}{Conf}

\newcommand\blfootnote[1]{%
	\begingroup
	\renewcommand\thefootnote{}\footnote{#1}%
	\addtocounter{footnote}{-1}%
	\endgroup
}
\title{A simple proof of the Bourgade-Sodin-Zeitouni theorem on the excess of the sine-$\beta$ process for non-integer $\sqrt{2/\beta}>1$}
\author{Sergei M. Gorbunov${}^*$
\blfootnote{Higher School of Modern Mathematics; Moscow Institute of Physics and Technology, Dolgoprudny, Moscow Region, Russia}
\blfootnote{Institute for System Programming of the Russian Academy of Sciences, Moscow, Russia}
\blfootnote{Steklov Mathematical Institute of Russian Academy of Sciences, Moscow, Russia}}

\date{}
\begin{document}
\maketitle
\begin{abstract}
We give an alternative proof of the Bourgade–Sodin–Zeitouni theorem on the excess of the sine-$\beta$ process for the Paley–Wiener space when $\sqrt{2/\beta}\in (1,+\infty)$ is not an integer. Let $k(\beta)=\lfloor\sqrt{2/\beta}\rfloor$. For almost every configuration, removing any $k(\beta)$ points leaves a minimal uniqueness set for the Paley-Wiener space. The same conclusion holds for each Bernstein space $\mathcal{B}^p$ of $L^p(\R)$-integrable functions with $m(\beta)<p<\infty$, where $m(\beta)=\beta(k(\beta)+1)-\sqrt{\beta^2(k(\beta)+1)^2-2\beta}$. 
\end{abstract}
\section{Introduction}
\subsection{Random configurations and minimal uniqueness sets}

For $\beta>0$, let $\P_\beta$ denote the sine-$\beta$ point process \cite{KillipStoiciu, ValkoViragA}. It is a Borel probability measure on the space of configurations $\Conf(\R)$ --- the space of locally finite subsets of the real line endowed with the vague topology. This process arises as the weak limit of the circular $\beta$-ensembles
\[
d\P_{\beta}^n(\theta_1, \ldots, \theta_n) = Z_{n, \beta}^{-1}\prod_{1\le l< m\le n}\abs{e^{i\theta_l}-e^{i\theta_m}}^\beta\prod_{k=1}^nd\theta_k, \quad \theta_j\in(-\pi, \pi)
\]
under the microscopic scaling $\theta_j\mapsto \frac{n}{2\pi}\theta_j$.

For $\beta=2$, the sine-$\beta$ process is determinantal with the sine kernel
\[
K_{\mathcal{S}}(x, y) = \frac{\sin\pi(x-y)}{\pi(x-y)}.
\]
By the Paley-Wiener Theorem, this is the integral kernel of the orthogonal projection from $L^2(\R)$ onto the \textit{Paley-Wiener space} $\mathcal{PW}$ --- the space of entire functions of exponential type at most $\pi$ whose restrictions to $\R$ belong to $L^2(\R)$.

A set $X\subset \R$ is a \textit{uniqueness set} for a function space $\mathcal{H}$ if the only function in $\mathcal{H}$ that vanishes on $X$ is the zero function. A uniqueness set $X$ is \textit{minimal} if every proper subset $Y\subset X$ admits a nonzero function $f\in\mathcal{H}$ vanishing on $Y$.

\begin{theorem}[Ghosh \cite{Ghosh}]
$\P_2$-almost every configuration $X\subset \R$ is a uniqueness set for the space $\mathcal{PW}$.
\end{theorem}

Bufetov, Qiu and Shamov \cite{BufetovQiuShamov} proved the Lyons-Peres conjecture. For a determinantal process induced by an orthogonal projection onto a reproducing kernel Hilbert space, almost every configuration is a uniqueness set for that space. The discrete analogue is due to Lyons \cite[Theorem 7.11]{Lyons}. A natural question is whether the random configuration is also minimal.

\begin{theorem}[Bufetov \cite{BufetovA}]\label{sine_excess}
For $\P_2$-almost every configuration $X\subset\R$ and every point $x\in X$, the subset $X\setminus\{x\}$ is a minimal uniqueness set for the Paley-Wiener space.
\end{theorem}

The \textit{excess} of a configuration relative to a function space is the number of points one must remove to obtain a minimal uniqueness set. Thus, relative to the Paley-Wiener space, the sine process has excess one almost surely. Recently, Bourgade, Sodin and Zeitouni announced \cite{BourgadeHuang} a beautiful result: $\P_\beta$-almost every configuration has excess
\[
k(\beta)=\left\lfloor\sqrt{\frac{2}{\beta}}\right\rfloor
\]
relative to the Paley-Wiener space \cite{BourgadeSodinZeitouni}.
We give an alternative proof when $\sqrt{2/\beta}$ is a non-integer greater than one.

\subsection{Statement of results}

For $0<p\le \infty$, let $\mathcal{B}^p$ denote the space of entire functions of exponential type at most $\pi$ whose restrictions to $\R$ belong to $L^p(\R)$, see Section \ref{sec:bernstein}.
The Paley-Wiener space corresponds to the case $p=2$. By Theorem \ref{nikolski_bernstein}, the inclusion
\[
\mathcal{B}^p\subset\mathcal{B}^q
\]
holds whenever $0<p\le q\le \infty$. Thus every uniqueness set for $\mathcal{B}^q$ is also a uniqueness set for $\mathcal{B}^p$. Set $$m(\beta) = \beta(k(\beta)+1)-\sqrt{\beta^2(k(\beta)+1)^2-2\beta}.$$

\begin{theorem}\label{excess}
Let $0<\beta\le 2$ and $m(\beta)<p<\infty$. If $\sqrt{2/\beta}\notin\N$, then for $\P_\beta$-almost every configuration $X\subset \R$ and for every subset $\{x_1, \ldots, x_{k(\beta)}\}\subset X$, the set
\[
X\setminus\{x_1, \ldots, x_{k(\beta)}\} 
\]
is a minimal uniqueness set for the space $\mathcal{B}^p$.
If $\sqrt{2/\beta}\in\N$, then either
\[
X\setminus\{x_1, \ldots, x_{k(\beta)}\} \quad \text{ or }\quad X\setminus\{x_1, \ldots, x_{k(\beta)-1}\}
\]
is a minimal uniqueness set for the space $\mathcal{B}^p$ for $\P_\beta$-almost every configuration $X$.
\end{theorem}

\subsection{Acknowledgements}
The author thanks Alexander Bufetov for introducing him to the approach developed in \cite{BufetovB} and for many helpful discussions.

The author is a winner of the BASIS Foundation Competition and is deeply grateful to the Jury
and the sponsors. The author is a winner of the all-Russia mathematical August Moebius contest of graduate and undergraduate student papers and thanks the jury and the board for the high praise of his work.

\section{Proof of Theorem \ref{excess}}
Let $X\subset \R$ be a configuration. Assume that the limit
\begin{equation}\label{eq:vp_sum}
\lim_{R\to\infty} \sum_{x\in X\cap[-R, R]}\frac{x}{1+x^2}
\end{equation}
exists and the sum
\begin{equation}\label{eq:abs_sum}
\sum_{x\in X}\frac{1}{1+x^2}
\end{equation}
converges. Then for each fixed $w\in\C\setminus\R$, the following Euler product
\[
G_X(z, w) = \lim_{R\to\infty}\prod_{x\in X\cap[-R, R]}\frac{x-z}{x-w}
\]
defines an entire function of $z$. If $X$ is a random configuration distributed according to the sine-$\beta$ process, the above Euler product is a nonzero random constant multiple of the Valk\'o-Vir\'ag secular function \cite{ValkoViragB}. Its zeros are precisely the points of the configuration $X$. The following criterion, extending the result of Bufetov~\cite{BufetovB}, determines the excess.
\begin{proposition}\label{excess_criterion}
Let $0<q\le p\le\infty$, $c>0$ and $j\ge 1$. Assume the following:
\begin{itemize}
\item the function $G_X(z, ic)$ has exponential type at most $\pi$;
\item the lower bound
\begin{equation}\label{eq:ex_lower}
\inf_{T\in\R}\abs{G_X(iT, ic)}e^{-\pi T}(1+\abs{T})^A >0
\end{equation}
exists for some $A>0$;
\item the conditions
\[
\frac{1}{1+\abs{\cdot}^j}G_X(ic+\cdot, ic) \notin L^p(\R) \quad \text{ and }\quad  \frac{1}{1+\abs{\cdot}^{j+1}}G_X(ic+\cdot, ic)\in L^q(\R)
\]
are satisfied.
\end{itemize}

Then for any subset $\{x_1, \ldots, x_j\}\subset X$, the set $X\setminus\{x_1, \ldots, x_j\}$ is a minimal uniqueness set for the space $\mathcal{B}^a$ with $a\in [q, p]$. If the third hypothesis is replaced by the assumptions
\[
\frac{1}{1+\abs{\cdot}^{j-1}}G_X(ic+\cdot, ic) \notin L^p(\R)\quad \text{ and }\quad  \frac{1}{1+\abs{\cdot}^{j+1}}G_X(ic+\cdot, ic)\in L^q(\R),
\]
then for any $a\in[q, p]$, either $$X\setminus\{x_1, \ldots, x_{j-1}\}, \quad \text{ or }\quad X\setminus\{x_1, \ldots, x_j\}$$ is a minimal uniqueness set for $\mathcal{B}^a$.
\end{proposition}

Here and subsequently we set $\cdot = t$ to be the integration variable over horizontal lines.

Proposition \ref{excess_criterion} is proved in Section \ref{sec:bernstein}. For the Paley-Wiener space, Proposition \ref{excess_criterion} has been established by Bufetov \cite{BufetovB}.

The sum \eqref{eq:abs_sum} converges absolutely $\P_\beta$-almost surely because the process has constant first correlation function and the function $(1+t^2)^{-1}$ is integrable.

Our main tool is the following consequence of Theorem \ref{clt}, proved in \cite{Gorbunov}.
\begin{theorem}\label{bounds}
Let $\beta\le 2$. Along some subsequence $R_k\to\infty$, the sum \eqref{eq:vp_sum} converges $\P_\beta$-almost surely. For $z, w, u, v\in\H$, define the expressions
\begin{equation}\label{eq:rho_def}
\tilde d(z, w; u, v) = \abs*{\frac{(\bar{w}-u)(\bar{z}-v)}{(\bar{w}-v)(\bar{z}-u)}}, \qquad d(z, w) =\tilde d(z, w, z, w) = \frac{\abs{z-\bar w}^2}{4\Im z \Im w}.
\end{equation} 
For $z_j, w_j\in\H$, $\gamma_j\in\R$, we have the upper bound
\[
\E_\beta \prod_{j=1}^k \abs{G_X(z_j, w_j)}^{\gamma_j} \le \exp\left(\pi\sum_{j=1}^k\gamma_j\Im[z_j-w_j]\right) \prod_{j=1}^k d(z_j, w_j)^{\frac{\gamma_p^2}{2\beta}}\prod_{1\le j<l\le k} \tilde d(z_j, w_j; z_l, w_l)^{\frac{\gamma_j\gamma_l}{\beta}}.
\]
Moreover, there are positive constants $e, a$ such that the inequalities
\[
a \le \frac{\E\abs{G_X(z, w)}^\gamma}{e^{\gamma \pi(\Im z-\Im w)}d(z, w)^{\frac{\gamma^2}{2\beta}}} \le 1
\]
hold whenever $\Im z\ge e$ and $\Im w\ge e$.
\end{theorem}
\begin{remark}
Observe that the absolute value of the Euler product $G_X$ is independent of the signs of $\Im z$ and $\Im w$, so the bounds still hold if the parameters are in the lower half-plane.
\end{remark}
Theorem \ref{bounds} is proved in Subsection \ref{sec:bounds}.

We apply Proposition \ref{excess_criterion} to the random configuration distributed according to the sine-$\beta$ process.
Theorem \ref{bounds} asserts the bound
\[
\E_\beta \abs{G_X(z, ic)}^\gamma \le e^{\gamma\pi(\Im z - c)}d(z, ic)^{\frac{\gamma^2}{2\beta}} \le Ce^{\gamma\pi(\Im z -c)}\left(\frac{1+\abs{z}^2}{\Im z}\right)^{\frac{\gamma^2}{2\beta}}
\]
for some constant $C$.
Therefore, for any $\varepsilon>0$, the integral
\[
\int_{\R^2}\frac{e^{-\gamma(\pi+\varepsilon)\abs{y}}}{1+\abs{t+iy}^{\gamma k}}\abs{G_X(t+iy, ic)}^\gamma dtdy
\]
is almost surely finite for sufficiently large $k$ and sufficiently small $\gamma$, because so is its expectation. Intersecting the full-measure events for positive rational $\varepsilon>0$ and using monotonicity in $\varepsilon$ yields the integrability for every $\varepsilon>0$ on one full-measure event. By Proposition \ref{bern_crit}, the function $G_{X\setminus\{x_1, \ldots, x_k\}}(z, ic)$ belongs to the Bernstein space $\mathcal{B}^\gamma$ and, therefore, has exponential type at most $\pi$. Therefore, so does the function $G_X(z, ic)$, since multiplying by a polynomial does not change the exponential type. This proves the first requirement of Proposition \ref{excess_criterion}. Moreover, by Proposition \ref{zeros_bern}, the limit $\eqref{eq:vp_sum}$ exists $\P_\beta$-almost surely.

The following lemma gives the required lower bound \eqref{eq:ex_lower}.
\begin{lemma}\label{vert_asymp}
For $\P_\beta$-almost every $X$, the asymptotic
\[
\log\abs{G_X(iT, ic)} - \pi(T-c) = o(\log T)
\]
holds for large $T$.
\end{lemma}
Lemma \ref{vert_asymp} is proved in Subsection \ref{sec:im_growth}.

It remains to establish two integrability claims. First, we prove that for $\beta\le 2$ and for any sufficiently small $\varepsilon>0$, the integrability condition
\begin{equation}\label{eq:upper_bound}
\frac{1}{1+\abs{t}^{k(\beta)+1}}G_X(t+ic, ic)\in L^{m(\beta)+\varepsilon}(\R)
\end{equation}
holds $\P_\beta$-almost surely. This is done in Subsection \ref{sec:integr}. Next, in Subsection \ref{sec:nonintegr} we establish that for all sufficiently large $p$ the non-integrability condition
\begin{equation}\label{eq:lower_bound}
\frac{1}{1+\abs{t}^{k(\beta)}}G_X(t+ic, ic)\notin L^{p}(\R)
\end{equation}
holds $\P_\beta$-almost surely if $\sqrt{2/\beta}\notin\N$ and the condition
\begin{equation}\label{eq:lower_bound_crit}
\frac{1}{1+\abs{t}^{k(\beta)-1}}G_X(t+ic, ic)\notin L^{p}(\R)
\end{equation}
is valid $\P_\beta$-almost surely in the critical cases $\sqrt{2/\beta}\in\N$. Theorem \ref{excess} then follows from Proposition \ref{excess_criterion}.

\subsection{Proof of the integrability \eqref{eq:upper_bound}}\label{sec:integr}
Factor the Euler product as follows
\[
G_X(t+ic, ic) = G_X(i\abs{t}+ic, ic)G_X(t+ic, i\abs{t}+ic).
\]
By Lemma \ref{vert_asymp}, the desired claim \eqref{eq:upper_bound} follows from the integrability
\[
\frac{e^{\pi \abs{t}}\abs{t}^{\delta}}{1+\abs{t}^{k(\beta)+1}}G_X(t+ic, ic+i\abs{t}) \in L^{m(\beta)+\varepsilon}(\R)
\]
for some $\delta>0$.

By Theorem \ref{bounds}, we have the bound
\begin{equation}\label{eq:norm_ineq}
\left[\frac{e^{\pi \abs{t}}\abs{t}^{\delta}}{1+\abs{t}^{k(\beta)+1}}\right]^\gamma\E_\beta\abs{G_X(t+ic, ic+i\abs{t})}^\gamma \le C(1+\abs{t})^{\frac{\gamma^2}{2\beta}+\gamma\delta -\gamma(1+k(\beta))}
\end{equation}
for some constant $C$ uniformly in $t\in\R$. Set $\gamma=m(\beta)+\varepsilon$. 
If the inequality
\begin{equation}\label{eq:pol_ineq}
\frac{\gamma^2}{2\beta}-\gamma(1+k(\beta))<-1
\end{equation}
holds, choose $\delta>0$ small enough that the exponent on the right-hand side of the bound \eqref{eq:norm_ineq} is still less than $-1$. Integration over $t$ then gives a finite expected $\gamma$-th power of the $L^\gamma(\R)$ norm or quasinorm, and hence almost-sure finiteness.

Since $m(\beta)$ is the smaller root of the polynomial $$h_\beta(v) = v^2/2\beta - v(1+k(\beta))+1,$$ we have $h_\beta(m(\beta)+\varepsilon)<0$ if $\varepsilon$ is sufficiently small, which is precisely the inequality \eqref{eq:pol_ineq}.

Inclusion \eqref{eq:upper_bound} is proved.
\subsection{The non-integrability condition \eqref{eq:lower_bound}}\label{sec:nonintegr}
Fix the constant $c$ to be sufficiently large, so that the last claim of Theorem \ref{bounds}, applied to $G_X(z, ic)$ with $\Im z\ge c$, holds. Normalize the Euler product and denote
\[
\widetilde{G}_X(z, w) = e^{\pi\Im(w-z)}G_X(z, w).
\]
Introduce the random integral
\[
I_R(p) = \int_1^2\abs{\widetilde{G}_X(tR+ic, iR)}^p dt.
\]
Suppose that \eqref{eq:lower_bound} fails with probability $1$. For any $\delta>0$, Lemma \ref{vert_asymp} yields the convergence
\begin{equation}\label{eq:mom_conv}
\lim_{R\to\infty}R^{1-\delta-k(\beta)p}I_R(p)=0
\end{equation}
for $\P_\beta$-almost every configuration $X$.

In what follows, we provide a lower bound on $I_R(p)$ that holds with positive probability and contradicts the convergence \eqref{eq:mom_conv}. By the H\"older inequality, lower order moments $0<\mu<\nu<p$ provide the bound from below
\begin{equation}\label{eq:holder_bound}
I_R(p) \ge \frac{I_R(\nu)^{\frac{p-\mu}{\nu-\mu}}}{I_R(\mu)^{\frac{p-\nu}{\nu-\mu}}}.
\end{equation}
For each $\nu\in\R$, the last claim of Theorem \ref{bounds} yields
\[
a \le \frac{\E_\beta I_R(\nu)}{R^{\frac{\nu^2}{2\beta}}} \le 1
\]
for some constant $a$. First, by Markov's inequality, for any $\epsilon>0$, there exists a sufficiently large $B$ that gives the upper bound
\[
\P_\beta\left(I_R(\mu) \le BR^{\frac{\mu^2}{2\beta}}\right)>1-\epsilon
\]
for all $R\ge 1$ with probability arbitrarily close to $1$.
Next, for $0<\nu<\sqrt{2\beta}$, there are constants $D, b>0$ such that the lower bound
\[
\P_\beta\left(I_R(\nu)\ge DR^{\frac{\nu^2}{2\beta}}\right) > b>0
\]
holds for all $R\ge 1$ with positive probability.
This follows from the Paley-Zygmund inequality and the following bound on fractional moments of $I_R(\nu)$.
\begin{lemma}\label{fractional_moments}
For any $0<\eta<1$ and any $0<\nu<\sqrt{\frac{2\beta}{1+\eta}}$, we have the bound
\[
\E_\beta (I_R(\nu))^{1+\eta}\le CR^{(1+\eta)\frac{\nu^2}{2\beta}}
\]
for some constant $C$.
\end{lemma}
Lemma \ref{fractional_moments} is established in Subsection \ref{sec:frac_mom}.

Since $\epsilon$ may be arbitrarily small, the two bounds above hold simultaneously with positive probability
\[
\P_\beta\left(I_R(\nu)\ge DR^{\frac{\nu^2}{2\beta}}, I_R(\mu)\le BR^{\frac{\mu^2}{2\beta}}\right) >b-2\epsilon>0.
\]
On this event, H\"older's inequality \eqref{eq:holder_bound} gives the lower bound
\[
I_R(p) \ge L R^{\frac{p(\nu+\mu)-\mu\nu}{2\beta}}
\]
for some positive constant $L$. If the numbers $\mu, \nu$ are chosen sufficiently close to $\sqrt{2\beta}$, this contradicts the convergence \eqref{eq:mom_conv}: as $\nu, \mu$ approach $\sqrt{2\beta}$, the exponent approaches
\[
\frac{p(\nu+\mu)-\mu\nu}{2\beta} \to p\sqrt{\frac{2}{\beta}}-1.
\]
Choose $\nu, \mu$ so that the difference between the left-hand side above and the limit is less than $\varepsilon'>0$. Then we have that
\[
R^{1-\delta - k(\beta)p}I_R(p) \ge \exp\left(\log R\left[p\left(\sqrt{\frac{2}{\beta}}-k(\beta)\right)-\varepsilon'-\delta\right]\right),
\]
where the exponent is positive if $\delta$ and $\varepsilon'$ are chosen sufficiently small, since $k(\beta)<\sqrt{2/\beta}$. This lower bound holds with probability bounded away from zero for large $R$, contradicting the convergence in probability implied by \eqref{eq:mom_conv} under the almost-sure integrability assumption.

The desired nonintegrability \eqref{eq:lower_bound} holds with positive probability. Translation ergodicity now implies that the event has probability one.

\begin{theorem}\label{ergodicity}
For $\beta\le 2$, sine-$\beta$ process is ergodic with respect to translations.
\end{theorem}
In section \ref{sec:ergodicity}, we present the argument of Assiotis and Suzuki \cite{AssiotisSuzuki} on how Theorem \ref{ergodicity} follows from the result of Dumaz and Malvy \cite{DumazMalvy}.

The condition \eqref{eq:lower_bound_crit} follows from the above argument with $k(\beta)$ replaced by $k(\beta)-1$.

\begin{remark}
Alternatively, one can use tail triviality proved by Assiotis, Suzuki \cite{AssiotisSuzuki}, and the rigidity proved by Chhaibi, Najnudel \cite{ChhaibiNajnudel} (see also \cite{DHLM}). Following Bufetov \cite{BufetovA, BufetovB}, this implies the ergodicity under compactly supported diffeomorphisms of the real line. By Proposition \ref{excess_criterion}, the excess of a configuration is invariant under this action.
\end{remark}

\section{The properties of the function $G_X$}
\subsection{Bounds on the moments: proof of Theorem \ref{bounds}}\label{sec:bounds}

Assume that $f$ is a compactly supported smooth function. Denote
\[
S_f(X) = \sum_{x\in X}f(x), \quad X\in\Conf(\R),
\]
and set
\[
\overline{S}_f(X) = S_f(X) - \E_\beta S_f.
\]
Introduce the $p$-order Sobolev scalar product
\[
\langle f, g\rangle_{H^p(\R)} = \int_\R\abs{\omega}^{2p}\hat{f}(\omega)\overline{\hat{g}(\omega)}d\omega, \quad \hat{f}(\omega) = \frac{1}{2\pi}\int_\R e^{-i\omega x}f(x)dx,
\]
which induces the respective $p$-order Sobolev semi-norm. We denote the latter by $\normHR{p}{f}$.
If $\beta\le 2$, we have \cite[Proposition 3.2]{Gorbunov} the inequality
\[
\E_\beta\abs{\overline{S}_f}^2 \le C\normHR{1/2}{f}^2
\]
holds for some constant $C$. By extension by continuity, this allows us to define the random variable $\overline{S}_f$ as the $L^2$-limit, if the semi-norm $\normHR{1/2}{f}$ is finite.
\begin{theorem}\label{clt}
Let $\beta\le 2$, and $f$ be a real-valued function with the semi-norm $\normHR{1/2}{f}$ being finite. Then we have the bound
\[
\E_\beta e^{\overline{S}_f} \le e^{\frac{1}{\beta}\normHR{1/2}{f}^2}.
\]
If the seminorm $\normHR{1}{f}$ is finite, we have the estimate
\[
\E_\beta e^{\overline{S}_f}  =  e^{\frac{1}{\beta}\normHR{1/2}{f}^2}\left[1+O\left(\normHR{1}{f}\right)\right].
\]
\end{theorem}

We apply Theorem \ref{clt} to the function
\[
f_{z, w}(x) = \log\abs{z-x}-\log\abs{w-x}.
\]
First, restrict this function to the interval $[-R, R]$. Extended by zero to $\R$, the latter is not $1/2$-Sobolev regular, so we also introduce its continuous approximation $f_{z, w}^R(x)$, which is equal to $f_{z, w}$ on the interval $[-R, R]$ and is linearly interpolated to zero on the intervals $R\le \abs{x}\le R+1$.

Since the first correlation function of the sine-$\beta$ process is equal to $1$, the limit
\[
\lim_{R\to\infty} \E_\beta\abs*{\sum_{x\in X\cap [-R, R]}f_{z, w}(x) - \sum_{x\in X}f_{z, w}^R(x)} \le \lim_{R\to\infty} 2\abs{f_{z, w}(R)} = 0
\]
vanishes. On the other hand, it is clear that the limits
\[
\lim_{R\to\infty}\norm{f_{z, w}^R-f_{z, w}}_{L^2(\R)} = \lim_{R\to\infty}\norm*{\frac{d}{dx}\left(f_{z, w}^R(x)-f_{z, w}(x)\right)}_{L^2(\R)}=0
\]
vanish. The inequality $\normHR{1/2}{f}^2\le \norm{f}_{L^2(\R)}^2+ \norm{f'}_{L^2(\R)}^2$ implies the convergence $f_{z, w}^R\to f_{z, w}$ in the Sobolev half-order seminorm and, therefore, the convergence
\[
\overline{S}_{f_{z, w}\I_{[-R, R]}} \overset{L^1(\Conf(\R), \P_\beta)}{\to} \overline{S}_{f_{z, w}},\quad \text{ as }R\to\infty.
\]
Finally, this convergence is almost sure along some sequence $R_k$. Now Theorem \ref{clt} is applicable to this almost sure limit.

Theorem \ref{bounds} then follows from Theorem \ref{clt} after we calculate the limit of the expectation and the Sobolev semi-norms. Below we establish the equalities
\begin{equation}\label{eq:average}
\lim_{R\to\infty}\E_\beta \sum_{x\in X\cap [-R, R]}(\log\abs{z-x}-\log\abs{w-x}) = \pi\Im(z-w),
\end{equation}
and
\begin{equation}\label{eq:variance}
\langle f_{z, w}, f_{u, v}\rangle_{H^{1/2}(\R)} = \frac{1}{2}\log\abs*{\frac{(\bar{w}-u)(\bar{z}-v)}{(\bar{w}-v)(\bar{z}-u)}}.
\end{equation}
The latter yields
\[
\exp\left(\frac{1}{\beta}\normHR{1/2}{\sum_{j=1}^k \gamma_jf_{z_j, w_j}}^2\right) = \prod_{1\le j, l\le k} d(z_j, w_j; z_l, w_l)^{\frac{\gamma_j\gamma_l}{2\beta}}.
\]
Setting $k=1$ and using that the first-order semi-norm is equal to the $L^2(\R)$-norm of the derivative
\[
\normHR{1}{f_{z, w}}^2 = \frac{1}{2\pi}\norm{f_{z, w}'}_{L^2(\R)}^2,
\]
we obtain the last assertion of Theorem \ref{bounds}: the above expression converges to zero uniformly in the real parts as the imaginary parts converge to infinity. 

\textbf{Formula for the expectation}

Express the limit of the integral as follows:
\begin{multline*}
\lim_{R\to\infty}\E_\beta \sum_{x\in X\cap [-R, R]}(\log\abs{z-x}-\log\abs{w-x}) = \lim_{R\to\infty}\int_{-R}^R(\log\abs{z-x}-\log\abs{w-x})dx\\
= \lim_{R\to\infty} \Re \int_{-R}^R\frac{z-w}{w-x}dx + \lim_{R\to\infty}\Re \int_{-R}^R \int_w^z \frac{w-u}{(w-x)(u-x)}dudx.
\end{multline*}
The second integral converges absolutely, so the limit may be omitted and the order of integration may be changed. For $w, u\in \H$ the integral over $x$ in the second term vanishes, while the first integral is equal to $-(z-w)i\pi$. Taking the real part, we deduce the equality \eqref{eq:average}.

\textbf{Formula for the variance}

We proceed to the identity \eqref{eq:variance}. To compute the Fourier transform of $f_{z, w}$, we differentiate
\[
f'_{z, w}(x) = \Re \left(\frac{1}{x-z}-\frac{1}{x-w}\right) = H_z(x)-H_w(x), \quad H_w(x) = \frac{x-\Re w}{(x-\Re w)^2+(\Im w)^2},
\]
and use the formula
\[
\widehat{H_w}(\omega) = -\frac{i}{2}\operatorname{sgn} \omega e^{-\abs{\omega}\Im w-i\omega\Re w}
\]
to deduce the expression
\[
\widehat{f_{z, w}}(\omega) = \frac{1}{2\omega}\operatorname{sgn}\omega\left[ e^{-\abs{\omega}\Im w-i\omega\Re w}-e^{-\abs{\omega}\Im z-i\omega \Re z}\right]
\]
for the Fourier transform of $f_{z, w}$. We obtain the following formula
\[
\langle f_{z, w}, f_{u, v}\rangle_{H^{1/2}(\R)} = \frac{1}{2}\Re\left[\int_0^\infty \frac{(e^{-i\bar{w}\omega} - e^{-i\bar{z}\omega})(e^{iv\omega}-e^{iu\omega})}{\omega}d\omega \right]
\]
for the scalar product. Equality \eqref{eq:variance} follows by substitution of the analytic continuation of the Frullani formula
\[
\int_0^\infty \frac{e^{-ut}-e^{-vt}}{t}dt = \log \left(\frac{v}{u}\right), \quad\text{ for }u, v \text{ such that } \Re v, \Re u>0,
\]
applied to the function $e^{-ax}$ for positive $a$.

\subsection{Growth on the imaginary axis: proof of Lemma \ref{vert_asymp}}\label{sec:im_growth}

Fix $\varepsilon > 0$. We show that the derivative of the function
\[
\widetilde{G}_X^\gamma(iT, ic)\, T^{-\varepsilon\gamma} 
= \Bigl|G_X(iT, ic)\, e^{-\pi(T-c)}\, T^{-\varepsilon}\Bigr|^\gamma
\]
with respect to $T$ is $\mathbb{P}_\beta$-almost surely absolutely integrable in the $T$-variable, provided $\gamma > 0$ is sufficiently small. Consequently, for any $\varepsilon > 0$, the function $\widetilde{G}_X(iT, ic)\, T^{-\varepsilon}$ is $\mathbb{P}_\beta$-almost surely bounded for $T \in \mathbb{R}$. Thus, this boundedness holds with probability one for all rational $\varepsilon > 0$. For such $\varepsilon$, we obtain  
\[
\limsup_{T\to\infty} \frac{\log \abs{\widetilde{G}_X(iT, ic)}}{\log T} \le \varepsilon.
\]
Applying the same argument to $\widetilde{G}_X^{-1}$ yields  
\[
\liminf_{T\to\infty} \frac{\log \abs{\widetilde{G}_X(iT, ic)}}{\log T} \ge -\varepsilon,
\]
which establishes the desired claim.

The derivative in question is  
\begin{equation}\label{eq:derivative_int}
\gamma\abs{\widetilde{G}_X(iT, ic)}^\gamma T^{-\gamma\varepsilon}\left\{\sum_{x\in X}\frac{T}{x^2+T^2}-\E_\beta\sum_{x\in X}\frac{T}{x^2+T^2} - \frac{\varepsilon}{T}\right\}.
\end{equation}
By Theorem~\ref{bounds}, the expectation of the last term satisfies  
\[
\gamma \varepsilon\, T^{-1-\varepsilon\gamma}\, 
\mathbb{E}_\beta \bigl|\widetilde{G}_X(iT, ic)\bigr|^\gamma 
\le C\, |T|^{\frac{\gamma^2}{2\beta} - 1 - \varepsilon\gamma}
\]
for some constant $C$. For sufficiently small $\gamma > 0$, the exponent is strictly less than $-1$, ensuring that the integral over $T$ converges.

Applying the Cauchy---Bunyakovsky---Schwarz inequality to the first two terms and invoking Theorem~\ref{bounds} gives  
\[
\E_\beta \abs{\widetilde{G}_X(iT, ic)}^\gamma T^{-\gamma\varepsilon}\abs*{\sum_{x\in X}\frac{T}{x^2+T^2}-\E_\beta\sum_{x\in X}\frac{T}{x^2+T^2}}\le \tilde C\abs{T}^{\frac{\gamma^2}{\beta}-\varepsilon\gamma-1}\normHR{1/2}{(\cdot/T)^2+1)^{-1}}.
\]
The scaling property of the half-order Sobolev seminorm
\[
\bigl\|f(\cdot/T)\bigr\|_{\dot H^{1/2}} = \|f\|_{\dot H^{1/2}}
\]
implies that the integral over $T$ has finite expectation if $\gamma$ is sufficiently small. Hence, the derivative \eqref{eq:derivative_int} is $\mathbb{P}_\beta$-almost surely integrable. Lemma~\ref{vert_asymp} follows.

\subsection{Fractional moments: proof of Lemma \ref{fractional_moments}}\label{sec:frac_mom}
We begin with an elementary integral inequality.

\begin{lemma}\label{frac_int_ineq}
Let $0<\eta<1$, let $h$ be a bounded, positive Borel function on $[1, 2]$. Then
\begin{equation}\label{eq:frac_int}
\left(\int_1^2 h(t)dt\right)^{1+\eta} \le \frac{1}{\log 2}\int_1^2\int_0^1 h(t)\left(\int\limits_{\substack{t'\in [1, 2] \\ u\le \abs{t'-t}\le 2u}}h(t')dt'\right)^\eta\frac{du}{u}dt.
\end{equation}
\end{lemma}
\begin{proof}
Without loss of generality, normalize $\int_1^2 h(t)dt = 1$, so the left-hand side of \eqref{eq:frac_int} equals $1$. The $t'$-integral on the right-hand side of \eqref{eq:frac_int} is bounded by the total integral of $h$, which is $1$.
Hence, using $b \le b^\eta$ for $0 < b \le 1$, we obtain
\begin{multline*}
\int_1^2\int_0^1 h(t)\left(\int\limits_{\substack{t'\in [1, 2] \\ u\le \abs{t-t'}\le 2u}}h(t')dt'\right)^\eta\frac{du}{u}dt \ge \int_1^2\int_0^1 h(t)\int\limits_{\substack{t'\in [1, 2] \\ u\le \abs{t-t'}\le 2u}}h(t')dt'\frac{du}{u}dt \\= \int_{[1, 2]^2}h(t)h(t')\int_{\abs{t-t'}/2}^{\abs{t-t'}}\frac{du}{u}dtdt' = \log 2.
\end{multline*}
\end{proof}

 We proceed to the proof of Lemma \ref{fractional_moments}.
 
Apply Lemma \ref{frac_int_ineq} to the fractional moment of $\widetilde{G}_X(Rt+ic, iR)$. For convenience, set
\[
z_R(t) = Rt + ic \in \H.
\]
For any $w\in \H$, we have
\[
\widetilde{G}_X(z_R(\cdot), iR) = \widetilde{G}_X(z_R(\cdot), w)\widetilde{G}_X(w, iR).
\]
Substitute $h(t) = \abs{\tilde{G}_X(z_R(t), iR)}^\nu$ into Lemma \ref{frac_int_ineq} together with the above factorization to obtain
\begin{equation*}
I_R^{1+\eta}(\nu) \le \frac{1}{\log 2}\int_1^2\int_0^1\abs{\widetilde{G}_X(w, iR)}^{(1+\eta)\nu}\abs{\widetilde{G}_X(z_R(t), w)}^\nu\times\left(\int\limits_{\substack{t'\in [1, 2]\\ u\le \abs{t-t'}\le 2u}}\abs{\widetilde{G}_X(z_R(t'), w)}^\nu dt'\right)^\eta \frac{du}{u}dt.
\end{equation*}
Observe that the point $w\in\H$ may depend on $t, u, R$, but must not depend on the $t'$ variable. Its choice will be specified later. 

Next, we use the H\"older inequality in the form
\[
\E AB^\eta \le (\E A)^{1-\eta}(\E AB)^\eta,
\]
to bound the average $\E_\beta I_R^{1+\eta}(\nu)$. Set
\begin{equation*}
A = \abs{\widetilde{G}_X(w, iR)}^{(1+\eta)\nu}\abs{\widetilde{G}_X(z_R(t), w)}^\nu\quad\text{ and }\quad
B = \int\limits_{\substack{t'\in [1, 2]\\ u\le \abs{t-t'}\le 2u}}\abs{\widetilde{G}_X(z_R(t'), w)}^\nu dt'.
\end{equation*}

Denote $a = \frac{\nu^2}{2\beta}$. Theorem \ref{bounds} yields bounds on the expectations of $A$ and $AB$:
\[
\E_\beta A \le E_A = d(z_R(t), w)^a d(w, iR)^{(1+\eta)^2a}\tilde d(z_R(t), w, w, iR)^{2(1+\eta)a},
\]
\[
\frac{\E_\beta AB}{E_A} \le E_B = \int\limits_{\substack{t'\in [1, 2]\\ u\le \abs{t-t'}\le 2u}} d(z_R(t'), w)^a\tilde d(z_R(t'), w, w, iR)^{2(1+\eta)a}\tilde d(z_R(t'), w, z_R(t), w)^{2a}dt'.
\]

First, observe the inequality
\[
\tilde d(z_1, z_2, z_2, z_3) = \tilde d(z_2, z_3, z_1, z_2) \le 10^{146}, \quad \text{ for }z_1, z_2, z_3\in\H.
\]
Next, if we set $w=z_R(t)+iuR$, we deduce the bounds
\[
d(z_R(t), w) + d(z_R(t'), w) \le C(1+uR), \quad \tilde{d}(z_R(t'), w, z_R(t), w) \le C, \quad d(w, iR) \le C \frac{R}{1+uR}
\]
for some positive constant $C$. We substitute the above bounds to obtain the inequality
\[
E_A \le C_A(1+uR)^a\left[\frac{R}{1+uR}\right]^{(1+\eta)^2a}, \quad E_B \le C_Bu(1+uR)^a
\]
for some constants $C_A, C_B$.

Substituting these bounds and integrating over $t$ yields
\begin{multline*}
\E_\beta I_R(\nu)^{1+\eta} \le \int_0^1 E_A E_B^\eta\frac{du}{u} \le  C_A C_B^{\eta}\int_0^1(1+uR)^{(1+\eta)a} \left[\frac{R}{1+uR}\right]^{(1+\eta)^2a}u^{\eta-1}du \\=C_A C_B^{\eta} R^{(1+\eta)a}\int_0^1\frac{u^{\eta-1}}{(u+R^{-1})^{\eta(1+\eta)a}}du
\le C_A C_B^{\eta}R^{(1+\eta)a}\int_0^1u^{\eta[1-(1+\eta)a]-1}du.
\end{multline*}
The conditions $\eta>0$ and $(1+\eta)a<1$, as assumed in the lemma, ensure that the final integral converges. This completes the proof of Lemma \ref{fractional_moments}.

\section{Ergodicity of the sine-$\beta$ process: proof of Theorem \ref{ergodicity}}\label{sec:ergodicity}

Following Assiotis, Suzuki \cite{AssiotisSuzuki}, Theorem \ref{ergodicity} follows from the result of Dumaz and Malvy \cite{DumazMalvy}. For completeness, we give their argument below.

For a configuration $X$ and a bounded Borel subset $B\subset \R$, set $\#_B(X)=\abs{X\cap B}$ to be the number of points. Let $\#_B$ be the respective random variable, with $X$ distributed according to the sine-$\beta$ process. Let $T_RB$ be the translation of the subset $B\subset \R$ on $R\in\R$.
\begin{theorem}\label{mixing_corr}
Let $I_1, \ldots, I_k$ and $J_1, \ldots, J_l$ be tuples of disjoint intervals. Then we have that
\[
\lim_{R\to\infty} \E_\beta \prod_{p=1}^k\#_{I_p}\prod_{q=1}^l\#_{T_RJ_q}=\E_\beta \prod_{p=1}^k\#_{I_p}\E_\beta \prod_{q=1}^l\#_{J_q}.
\]
That is, the sine-$\beta$ process is mixing on products of numbers of points in disjoint intervals.
\end{theorem}

In order to establish the mixing property on $L^2(\P_\beta)$, it is sufficient to consider that property on the sets
\[
N_B = \{X\in \Conf(\R): X\cap B=\emptyset\}, \quad \text{ for a finite union of intervals }B\subset \R.
\]
Indeed, these subsets form a $\pi$-system. Let us show that they generate the $\sigma$-algebra of measurable subsets of $\Conf(\R)$. Recall that the latter is the minimal $\sigma$-algebra with respect to which the functions $\#_I$, with $I\subset \R$ ranging among finite intervals, are measurable. On simple configurations these indicators are the limits
\[
\#_I = \lim_{n\to\infty}\sum_{j=1}^n(1-\I_{N_{I_j}}), \quad \text{ where } I=\cup_{j=1}^nI_j
\]
is the partition of $I$ into intervals of length $o(1)$ for large $n$. We conclude that the $\sigma$-algebra generated by the sets $N_B$ coincides with the Borel $\sigma$-algebra of the vague topology on simple configurations. It remains to use the following lemma.
\begin{lemma}
Let $\mathcal{A}$ be a $\pi$-system of subsets of $\Omega$ containing the full set $\Omega$. Set $\mathcal{C}= \sigma(\mathcal{A})$, and let $\P$ be a probability measure on the $\sigma$-algebra $\mathcal{C}$. Then the linear span of the indicators $\{\I_A\}_{A\in\mathcal{A}}$ is dense in $L^2(\P, \mathcal{C})$.
\end{lemma}
\begin{proof}
Indeed, consider a family of measurable subsets, on which the integral of a function $f\in L^2(\P, \mathcal{C})$ vanishes
\[
\int_Bfd\P=0.
\]
The family of such subsets $B$ is closed under countable disjoint unions, and, if the above equality holds for $B=\Omega$, forms a $\lambda$-system. Therefore, if $f$ is orthogonal to the indicator functions $\{\I_A\}_{A\in\mathcal{A}}$, then it is orthogonal to the indicator functions of subsets of the corresponding $\lambda$-system. The claim now follows from Dynkin's lemma.
\end{proof}

Let us now establish the mixing property on the subsets $N_B$. For every configuration, the indicator function of the set $N_I$ is the pointwise limit
\begin{equation}\label{eq:appr}
\I_{N_I}(X) = \lim_{n\to\infty}\prod_{j=1}^n(1-\#_{I_j}), \quad \text{ where } I=\cup_{j=1}^nI_j
\end{equation}
is the partition of $I$ into intervals of length $o(1)$ for large $n$. Since the sine-$\beta$ process is simple, this is also $\P_\beta$-almost sure limit. Further, we have the bound
\[
\abs*{\prod_{j=1}^n(1-\#_{I_j})-N_I} \le 5e^{2\#_I},
\]
which is uniform in $n$.
The function on the right-hand side may be bounded by an exponential of an additive functional
\[
S_f(X)=\sum_{x\in X}f(x),
\]
with $f$ smooth. The expectation $\E e^{2 S_f}$ is finite by Theorem \ref{clt}, so the approximation \eqref{eq:appr} holds in $L^2(\P_\beta)$ and, by translation invariance of $\P_\beta$, is uniform under translations of the subsets entering the relation. Theorem \ref{mixing_corr} implies the equality
\[
\lim_{R\to\infty}\E_\beta \prod_{p=1}^n(1-\#_{I_p})\prod_{q=1}^n(1-\#_{T_RJ_q}) = \E_\beta \prod_{p=1}^n(1-\#_{I_p})\E_\beta\prod_{q=1}^n(1-\#_{J_q}),
\]
which yields the assertion of Theorem \ref{ergodicity}.

\begin{remark}
This proof relies on the simplicity of the sine-$\beta$ process: no point is occupied by two particles almost surely. For a general point process, simplicity follows from the existence of the first two correlation functions with respect to Lebesgue measure, see \cite[Proposition 5.4.VI]{DaleyVerejones}. The sine-$\beta$ process admits all correlation functions, see Assiotis, Najnudel \cite{AssiotisNajnudel}. The existence of the second correlation function also follows from the work of Qu, Valk\'o \cite{QuValko}. 

An elementary way to deduce the existence of the correlation functions is to use the characterization of $\P_\beta$ as the weak microscopic limit of the circular $\beta$-ensembles, and apply the Aomoto variation of the Selberg integral.
\end{remark}

\section{Bernstein spaces}\label{sec:bernstein}
Recall that an entire function $f$ has exponential type at most $\pi$ if for any $\varepsilon>0$, there exists a constant $C_\varepsilon>0$ such that the bound
\[
\abs{f(z)}\le C_{\varepsilon}e^{(\pi+\varepsilon)\abs{z}}
\]
holds. For $0<p\le \infty$, the Bernstein space $\mathcal{B}^p$ is the space of entire functions of exponential type at most $\pi$ whose restrictions to $\R$ belong to $L^p(\R)$. By the result of Plancherel and P\'olya \cite[Equation (59) on page 120]{PlancherelPolya}, the integrability assumption may be equivalently imposed on any horizontal line: complex shifts preserve the space $\mathcal{B}^p$.
\begin{lemma}\label{PlancherelPolya}
Assume that an entire function $f$ has finite exponential type. Then for any $p>0$ and any $c_1, c_2\in\R$ the conditions
\[
f(\cdot+ic_1)\in L^p(\R) \quad \text{ and }\quad f(\cdot+ic_2)\in L^p(\R)
\]
are equivalent.
\end{lemma}

The Bernstein spaces are related by the Nikolskii-Bernstein inequalities.
\begin{theorem}\label{nikolski_bernstein}
For $0< p\le q\le \infty$, there exists a constant $C$ such that the inequality
\[
\norm{f}_{L^q(\R)} \le C\norm{f}_{L^p(\R)}
\]
holds for any function $f$ of exponential type $\pi$. In particular, we have the embeddings
\[
\mathcal{B}^p\subset\mathcal{B}^q.
\]
\end{theorem}
See Nessel and Wilmes \cite[Theorem 5]{NesselWilmes}.

The Nikolskii-Bernstein inequalities imply that every $f\in\mathcal{B}^p$ belongs to the Cartwright class: it has finite exponential type and the integral
\[
\int_\R \frac{\log_+\abs{f(x)}}{1+x^2}dx<+\infty
\]
is finite since $f\in L^\infty(\R)$. As states the Cartwright-Levinson Theorem, this yields the convergence of inverse zeros, see Levin \cite[Lecture 17.2, Theorem 1]{Levin}.
\begin{proposition}\label{zeros_bern}
For a function $f\in\mathcal{B}^p$, let $\rho$ range over its zeros. Then the limit
\[
\lim_{R\to\infty} \sum_{\rho: \abs{\rho}\le R, \rho\ne 0}\frac{1}{\rho}
\]
exists.
\end{proposition}

Let us establish a sufficient condition for a function to belong to a Bernstein space.
\begin{proposition}\label{bern_crit}
Let $0<p<\infty$. Let $f$ be an entire function. Assume that for any $\varepsilon>0$, the integral
\[
I_\varepsilon(f)=\int_{\R^2}e^{-p(\pi+\varepsilon)\abs{y}}\abs{f(x+iy)}^pdxdy
\]
is finite. Then $f$ belongs to $\mathcal{B}^p$.
\end{proposition}
\begin{proof}
Let us first show that the function has exponential type $\pi$. The function $\abs{f(z)}^p$ is subharmonic, so it is bounded by the average
\[
\abs{f(z)}^p \le \frac{1}{2\pi}\int_0^{2\pi}\abs{f(z+re^{i\theta})}^pd\theta
\]
for any real $r>0$, and, therefore, by the two-dimensional integral
\[
\abs{f(z)}^p \le \frac{1}{2\pi}\int_1^2\int_0^{2\pi}\abs{f(z+re^{i\theta})}^pd\theta rdr \le Ce^{p(\pi+\varepsilon)\abs{\Im z}}I_\varepsilon(f)
\]
for any $\varepsilon>0$ and some constant $C$. This proves the exponential type $\pi$.

The hypothesis of the proposition implies that for some $t\in \R$, we have $f(\cdot+it)\in L^p(\R)$. Given finite exponential type, Lemma \ref{PlancherelPolya} implies that it belongs to $\mathcal{B}^p$.
\end{proof}

We proceed to the proof of Proposition \ref{excess_criterion}.
\begin{proof}[Proof of Proposition \ref{excess_criterion}]
Once it is established that $X$ with $j$ points removed is a uniqueness set, its minimality follows from Lemma \ref{PlancherelPolya}. Let us show the former claim. Our argument follows Levin \cite[Lecture 18, Theorem 4]{Levin}.

Assume that there exists a function $g\in \mathcal{B}^a$, $a\in[q, p]$ vanishing on $X\setminus \{x_1, \ldots, x_j\}$. Consider the entire fraction
\[
u(z) = \frac{g(z)}{G_{X\setminus\{x_1, \ldots, x_j\}}(z)}.
\]
The claim follows once we show that the function $u$ has exponential type zero. Indeed, the lower bound \eqref{eq:ex_lower} on $G_X(iT)$ and the upper exponential bound on $g(iT)$ imply that $u(iT)$ has at most polynomial growth for large $T$. By the Phragm\'en-Lindel\"of principle, these two facts yield that $u$ is a polynomial \cite[Lecture 6.1, Problem 3]{Levin}. However, the requirement of the proposition and Lemma \ref{PlancherelPolya} imply that the function $G_{X\setminus\{x_1, \ldots, x_j\}}$ does not belong to $L^p(\R)$. By Theorem \ref{nikolski_bernstein}, neither does it belong to $L^a(\R)$. On the other hand, if multiplied by the polynomial $u$, it belongs to $L^a(\R)$. This leaves us with $u=0$ and, hence, $g=0$.

Below, we show that the function $u$ belongs to the Cartwright class. To a Cartwright function $v$, the Krein Theorem \cite[Lemma 3 and Equation (20)]{Krein} applies: there exist constants $\mu_1, \mu_2>0$, such that for any $\varepsilon>0$, the inequalities
\begin{align*}
\frac{\log\abs{v\left(re^{i\theta}\right)}}{r}\le \mu_1\sin\theta + \varepsilon,\\
\frac{\log\abs{v\left(re^{-i\theta}\right)}}{r}\le \mu_2\sin\theta + \varepsilon
\end{align*}
hold for all $\theta\in[0, \pi]$ and all sufficiently large $r$. Further, the indicator function of $v$ is equal to $\mu_1\sin\theta$ or $\mu_2\sin\theta$ in the upper and lower half-plane respectively, see Levin \cite[Lecture 16.1, Formula (7)]{Levin}. The constants $\mu_1, \mu_2$ corresponding to $u$ are zero due to the polynomial bound on the imaginary axis. This implies that the exponential type of $u$ is zero.

To show that $u$ belongs to the Cartwright class, we use Krein's characterisation of that class \cite[Lecture 16.1, Theorem 1]{Levin}: a function $v$ is Cartwright if and only if $\log\abs{v(z)}$ admits positive harmonic majorants in the lower/upper half-plane. Both the function $g$ and the function $G_{X\setminus\{x_1, \ldots, x_j\}}$ belong to the Cartwright class, since $g\in \mathcal{B}^a$ and $G_{X\setminus\{x_1, \ldots, x_{j+1}\}}\in \mathcal{B}^q$. Let $A$ be the positive harmonic majorant of $\log\abs{g}$, and $B$ be the positive harmonic majorant of $\log\abs{G_{X\setminus\{x_1, \ldots, x_j\}}}$ in the upper/lower half-plane. Since the function $G_{X\setminus\{x_1, \ldots, x_j\}}$ has only real zeros, the inequality
\[
\log\abs{u(z)}\le A + B - \log\abs*{G_{X\setminus\{x_1, \ldots, x_j\}}}
\]
provides a positive harmonic majorant for $\log\abs{u}$ in each half-plane.

For the last part of Proposition \ref{excess_criterion}, observe that the two possibilities
\[
G_{X\setminus\{x_1, \ldots, x_j\}}\in L^a(\R), \quad \text{ or }\quad G_{X\setminus\{x_1, \ldots, x_j\}}\notin L^a(\R)
\]
correspond to the set $X\setminus\{x_1, \ldots, x_{j-1}\}$ or the set $X\setminus\{x_1, \ldots, x_{j}\}$ being the minimal uniqueness set by the above argument.

Proposition \ref{excess_criterion} is proved.
\end{proof}


\begin{thebibliography}{99}

\bibitem{AssiotisNajnudel}
T.~Assiotis and J.~Najnudel,
\textit{Moments of C$\beta$E field partition function,
$\mathsf{Sine}_{\beta}$ correlations and stochastic zeta},
arXiv:2602.08739 (2026).

\bibitem{AssiotisSuzuki}
T.~Assiotis and K.~Suzuki,
\textit{Sine beta is tail trivial},
arXiv:2610.10849 (2026).

\bibitem{BourgadeHuang}
P.~Bourgade, J.~Huang,
\textit{Loop Equations Characterize Random Matrix Statistics},
arXiv:2607.07617 (2026).

\bibitem{BourgadeSodinZeitouni}
P.~Bourgade, M.~Sodin, and O.~Zeitouni,
\textit{The Paley--Wiener excess and a local law for the
$\mathrm{Sine}_{\beta}$ process},
private communication.

\bibitem{BufetovA}
A.~I.~Bufetov,
\textit{The sine-process has excess one},
arXiv:1912.13454 (2019).

\bibitem{BufetovB}
A.~I.~Bufetov,
\textit{Stochastic Euler products in the Paley--Wiener space},
arXiv:2611.01347 (2026).

\bibitem{BufetovQiuShamov}
A.~I.~Bufetov, Y.~Qiu, and A.~Shamov,
\textit{Kernels of conditional determinantal measures and the
Lyons--Peres completeness conjecture},
J. Eur. Math. Soc. \textbf{23} (2021), no.~5, 1477--1519.

\bibitem{ChhaibiNajnudel}
R.~Chhaibi and J.~Najnudel,
\textit{Rigidity of the $\mathrm{Sine}_{\beta}$ process},
Electron. Commun. Probab. \textbf{23} (2018), paper no.~94, 8~pp.

\bibitem{DaleyVerejones}
D.~J.~Daley and D.~Vere-Jones,
\textit{An Introduction to the Theory of Point Processes.
Vol.~I: Elementary Theory and Methods},
2nd ed., Springer, New York, 2003.

\bibitem{DHLM}
D.~Dereudre, A.~Hardy, T.~Lebl\'e, and M.~Ma\"ida,
\textit{DLR equations and rigidity for the sine-beta process},
Comm. Pure Appl. Math. \textbf{74} (2021), no.~1, 172--222.

\bibitem{DumazMalvy}
L.~Dumaz and M.~Malvy,
\textit{Long-range correlation of the $\mathrm{Sine}_{\beta}$ point process},
arXiv:2603.15289 (2026).

\bibitem{Ghosh}
S.~Ghosh,
\textit{Determinantal processes and completeness of random exponentials:
the critical case},
Probab. Theory Related Fields \textbf{163} (2015), nos.~3--4, 643--665.

\bibitem{Gorbunov}
S.~M.~Gorbunov,
\textit{Gessel-type expansion for the circular $\beta$-ensemble and
central limit theorem for the Sine-$\beta$ process for $\beta \leq 2$},
arXiv:2510.15172 (2025).

\bibitem{KillipStoiciu}
R.~Killip and M.~Stoiciu,
\textit{Eigenvalue statistics for CMV matrices: from Poisson to clock
via random matrix ensembles},
Duke Math. J. \textbf{146} (2009), no.~3, 361--399.

\bibitem{Krein}
M.~G.~Krein,
\textit{A contribution to the theory of entire functions of exponential type},
Izv. Akad. Nauk SSSR Ser. Mat. \textbf{11} (1947), no.~4, 309--326
(in Russian).

\bibitem{Levin}
B.~Ya.~Levin,
\textit{Lectures on Entire Functions},
in collaboration with Yu.~Lyubarskii, M.~Sodin, and V.~Tkachenko,
Translations of Mathematical Monographs, vol.~150,
American Mathematical Society, Providence, RI, 1996.

\bibitem{Lyons}
R.~Lyons,
\textit{Determinantal probability measures},
Publ. Math. Inst. Hautes \'Etudes Sci. \textbf{98} (2003), 167--212.

\bibitem{NesselWilmes}
R.~J.~Nessel and G.~Wilmes,
\textit{Nikolskii-type inequalities for trigonometric polynomials
and entire functions of exponential type},
J. Austral. Math. Soc. Ser.~A \textbf{25} (1978), no.~1, 7--18.

\bibitem{PlancherelPolya}
M.~Plancherel and G.~P\'olya,
\textit{Fonctions enti\`eres et int\'egrales de Fourier multiples. II},
Comment. Math. Helv. \textbf{10} (1937--1938), 110--163.

\bibitem{QuValko}
Y.~Qu and B.~Valk\'o,
\textit{On the pair correlation function of the $\mathrm{Sine}_{\beta}$ process},
arXiv:2509.15446 (2025).

\bibitem{ValkoViragA}
B.~Valk\'o and B.~Vir\'ag,
\textit{Operator limit of the circular beta ensemble},
Ann. Probab. \textbf{48} (2020), no.~3, 1286--1316.

\bibitem{ValkoViragB}
B.~Valk\'o and B.~Vir\'ag,
\textit{The many faces of the stochastic zeta function},
Geom. Funct. Anal. \textbf{32} (2022), no.~5, 1160--1231.

\end{thebibliography}
\end{document}